\input amstex
\documentstyle{grubb-p}

\def\d{d\!{}^{\!\text{\rm--}}\!}

\def\crp{\overline{\Bbb R}_+}

\def\rnp{{\Bbb R}^n_+}

\def\Zfrac{\tsize\frac1{\raise 1pt\hbox{$\scriptstyle z$}}}
\def\zfrac{\frac1{\raise 1pt\hbox{$\scriptscriptstyle z$}}}

\define\res{\operatorname{res}}
\define\tr{\operatorname{tr}}

\define\Tr{\operatorname{Tr}}
\define\TR{\operatorname{TR}}

\document

\topmatter
\title
{Analysis of Invariants Associated with Spectral Boundary
Problems for Elliptic Operators}
\endtitle
\author{Gerd Grubb}
\endauthor
\address
Copenhagen Univ. Math. Dept.,
Universitetsparken 5, DK-2100 Copenhagen, Denmark.
\endaddress
\email  grubb\@math.ku.dk \endemail
\subjclassyear{2000} \subjclass {58J42, 35S15, 58J32, 41A60}
\endsubjclass
\rightheadtext{Analysis of invariants}
\thanks
{This work was supported in part by The Danish Science Research Council,
SNF grant 21-02-0446}
\endthanks
\endtopmatter

\head Introduction \endhead
\medskip

Boundary conditions defined by pseudodifferential projections ---
also called spectral boundary conditions ---$$
\Pi u|_{\partial  X}=0\tag1
$$ for Dirac
operators $D$ were first 
used by Atiyah, Patodi and Singer in their important paper
\cite{APS1}, which also introduced the eta invariant. Such problems have
since then been studied in numerous other works 
(of which a section is found in the references to this survey).
The operator $(D_\Pi)^*D_\Pi 
$ is a realization of the Laplace operator $D^*D$ with complementing
projection boundary conditions $$
\Pi u|_{\partial  X}=0,\quad \Pi ^\perp
Du|_{\partial  X}=0.\tag 2
$$ More generally, one can consider a second order operator $P$ of
Laplace-type together with a boundary condition similar to (2); such
problems have 
been studied by the author in \cite{G6} (for the motivation in
physics, see the 
introduction there and Vassiliev \cite{V1, V2}).

In this survey paper we shall give an account of recent results concerning
some of the basic geometric invariants associated with such
operators. They are defined in analysis as 
coefficients in trace expansions for associated heat operators or
resolvents, expanded in powers and logarithms of the time variable
$t$ or spectral 
variable $\lambda $; they can also be defined from the pole structure
of meromorphic extensions of associated zeta and eta functions.

We shall consider two basic questions: 1) What happens to the
coefficients (or specific ones of them) when the boundary projection
is changed?
2) What happens to the coefficients when the 
interior operator is changed?

In both cases, we focus particularly on the coefficients of 
logarithmic terms
and the global coefficients ``behind'' them.

\example{Remark} The use of the terminology ``spectral boundary condition'',
earlier applied 
to boundary conditions defined by the positive eigenprojection of a
first-order elliptic differential operator, is motivated here by the fact that
{\it any} orthogonal 
pseudodifferential projection is the positive eigenprojection for a
first-order elliptic pseudodifferential operator (see Theorem 1.2 below). 
\endexample

\bigskip

\head 1. Pseudodifferential operators on closed manifolds \endhead

\medskip\subhead 1.1 Trace expansions \endsubhead

To prepare for the case of manifolds with boundary, we first recall
some results for closed manifolds.

Consider an $n$-dimensional compact $C^\infty $
manifold $X$ without boundary. 
Let $A$ be a classical (also called one-step polyhomogeneous)
pseudodifferential operator ($\psi $do) of
order $\nu \in\Bbb R$, acting on the sections of a $C^\infty $ vector
bundle $E$ over $X$ of dimension $n_1$. Let $P$ be a
classical elliptic $\psi $do of positive integer order $m$,
likewise acting in $E$ and such that the principal symbol has no
eigenvalues on $\Bbb R_-$. As shown by Grubb and Seeley in
\cite{GS1, Th\. 2.7} for (1.1), with the transition to the
essentially equivalent statements (1.2) and (1.3)
accounted for e.g\. in \cite{GS2}, a calculation in local coordinates
gives trace expansions
(we denote $\{0,1,2,\dots\}=\Bbb N$):$$
\gather
\Tr\bigl( A(P-\lambda )^{-N}\bigr)
\sim
\sum_{ j\in \Bbb N  } \tilde c_{ j}(-\lambda ) ^{\frac{\nu +n -j}m-N}+ 
\sum_{k\in \Bbb N}\bigl( \tilde c'_{ k}\log (-\lambda ) +\tilde c''_{
k}\bigr)(-\lambda ) ^{{ -k}-N},
\tag1.1\\
\Gamma (s)\Tr (AP^{-s})\sim
\sum_{ j\in \Bbb N } \frac{c_{j}}{s+\frac{j-\nu -n}m}-\frac{\Tr (A\Pi
_0(P))}s +  \sum_{k\in \Bbb N}\Bigl(\frac{
c'_{k}}{(s+k)^2}+\frac{c''_{k}}{s+ k}\Bigr) ,\tag1.2\\
\Tr  (Ae^{-tP})\sim
\sum_{ j\in \Bbb N } c_{j}t^{\frac{j-\nu -n}m}+ 
\sum_{k\in \Bbb N}\bigl({- c'_{k}}\log t+{c''_{k}}\bigr) t^{{k}}.\tag1.3
\endgather
$$ 
In (1.1), $N>\frac {\nu +n}m$ and $\lambda \to\infty $ on rays in 
an open subsector of $\Bbb
C$ containing $\Bbb R_-$. (1.2) 
means that $\Gamma (s)\Tr (AP^{-s})$,
defined in a standard way for $\operatorname{Re}s>\frac {\nu +n}m$, 
extends meromorphically to $\Bbb C$ with the pole structure 
indicated in the right hand side. Here $\Pi _0(P)$ is the orthogonal
projection onto
the nullspace $V_0(P)$ of $P$ (on which $P^{-s}$ is taken to be zero).
(1.3) holds when the resolvent for large $\lambda $ is defined for
$|\arg \lambda -\pi 
|\le\frac \pi 2+\varepsilon $, some $\varepsilon >0$; here $t\to 0+$
and the coefficients are {\it the same} as those in (1.2). 
There are universal nonzero proportionality factors linking the 
coefficients $\tilde c_j$ and $c_j$, $\tilde c'_k$ and $ c'_k$,
resp\. $\tilde c''_k$ and $c''_k$. 

A remarkable feature is the presence of the series over $k$ that
appears when $A$ and $P$ are not differential operators. 
Each of the coefficients $\tilde c_j$ and $\tilde c_k'$ comes from a specific
homogeneous term in the symbol of $A(P-\lambda)^{-N}$, whereas the
coefficients 
$\tilde c''_k$ depend on the full symbol. Thus the
coefficients $\tilde c_j$ and $\tilde c_k'$ depend each on a finite
set of homogeneous terms in the symbols of $A$ and $P$; we call such
coefficients `locally determined' (or `local'), while the $\tilde
c''_k$ are called `global'. When $\nu \notin
\Bbb Z$, the $\tilde c'_k$ vanish. When $\nu \in\Bbb Z$ and $(j-n-\nu
)/m$ is an integer $k\ge 0$, both $\tilde c_j$ and $\tilde c''_k$
contribute to the power $(-\lambda )^{-k-N}$; their sum is
independent of the choice of local coordinates, whereas the
splitting in $\tilde c_j$ and $\tilde c''_k$ depends in a well-defined way
on the symbol structure in local coordinates (see  \cite{GS1,
Th.\ 2.1} or \cite{G7, Th. 1.3}).

Particular attention has been paid to the first two coefficients
$\tilde c'_0$ and $\tilde c''_0$ in the series over $k$. Here $\tilde
c'_0$ is proportional 
to the {\it
noncommutative residue} of $A$, as introduced by Wodzicki \cite{W}
and Guillemin \cite{Gu}:$$
\operatorname{res} A=m\cdot \operatorname{Res}_{s=0}\Tr (AP^{-s})
=m\cdot  c'_0;\tag1.4
$$
it is an integral of the ($-n$)th order symbol of $A$.
As for $\tilde c''_0$, it is only the full coefficient of $(-\lambda
)^{-N}$,  called $C_0(A,P)$,
$$
C_0(A,P)=\tilde c''_0+\tilde c_{\nu +n} ,\tag1.5
$$
that is determined from the operators $A$ and $P$ when $\nu $ is
integer $\ge -n$; the
splitting in two terms depends on the choice of
local coordinates. For convenience, we define $\tilde c_{\nu +n }=0$
when 
$\nu \notin\Bbb Z$ or $\nu
<- n$. Then we have in local coordinates:$$
 \tilde c'_0=c'_0,\quad \tilde c''_0=c''_0,\quad \tilde c_{\nu
+n}=c_{\nu +n}.\tag1.6 
$$

The coefficient $C_0(A,P)$
equals the canonical trace $\TR A$ introduced by Kontsevich and Vishik
\cite{KV}  when $\nu <-n$ (then it is simply $\Tr A$),
when $\nu \notin \Bbb Z$ (cf\. also Lesch \cite{L}, \cite{G7}), and
when $A$ and $P$ have symbols of 
even-even parity  (see below) and 
the dimension of the manifold is odd. 
The cases of \cite{KV} are supplied by the case of even-odd symbols 
on manifolds of even dimension in
\cite{G7} in this volume, where we derive the results in all cases by
resolvent expansion techniques.

We say (as in \cite{G5}, \cite{G7}) that a classical $\psi $do $Q$ of order
$r\in\Bbb Z$ with symbol $q\sim 
\sum_{l\in\Bbb N}q_{r-l}(x,\xi )$ ($q_{r-l}$ homogeneous of degree
$r-l$ in $\xi $ for $|\xi |\ge 1$) in a local coordinate system has {\bf even-even} alternating
parity (in short: is even-even), when the 
symbol terms with even (resp\. odd) degree are even (resp\. odd) in $\xi $:
$$
\partial _x^\beta \partial _\xi ^\alpha 
q_{r-l}(x,-\xi )=(-1)^{r-l-|\alpha |}\partial _x^\beta \partial _\xi
^\alpha q_{r-l}(x,\xi ) \text{ for }|\xi |\ge 1, \text{ all }\alpha
,\beta .\tag1.7
$$ 
The operator (or symbol) is said to have {\bf even-odd} alternating
parity in the reversed 
situation where the
symbol terms with even (resp\. odd) degree are odd (resp\. even) in $\xi $:
$$
\partial _x^\beta \partial _\xi ^\alpha q_{r-l}(x,-\xi
)=(-1)^{r-l-1-|\alpha |}\partial _x^\beta \partial _\xi ^\alpha q_{r-l}(x,\xi ) \text{ for }|\xi |\ge 1.\tag1.8
$$ 
\cite{KV} calls the even-even
symbols ``odd-class'', studying them on odd-di\-men\-si\-o\-nal manifolds.

More generally,
it can be shown that
$$
C_0(A,P)-C_0(A,P')\text{ and }C_0([A,A'],P)\text{ are locally
determined} \tag1.9
$$
(we call $C_0(A,P)$ a
quasi-trace then). In \cite{G7} this is deduced from resolvent
expansions. Using more functional calculus, involving properties of complex
powers and logarithms of $P$, one can show residue formulas for the
expressions in (1.9); this was done for the first expression in Okikiolu
\cite{O} and \cite{KV}, and for  
both expressions in Melrose
and Nistor \cite{MN}, with extensions to the $b$-calculus. In the
latter work and in sequels to it, 
$C_0(A,P)$ is called a regularized trace of $A$ and denoted e.g.\
$\widehat{\Tr}A$, $\Tr_P(A)$. It is studied from a physics point of
view in Cardona, Ducourtioux, Magnot, Paycha \cite{CDMP},
\cite{CDP}, where it is called a weighted trace.  

A generalization to boundary value
problems (in the Boutet de Monvel calculus) has been worked out in
a joint work with Schrohe \cite{GSc}, on the basis of 
resolvent considerations (the complex powers are not inside the calculus). 

By division by $\Gamma (s)$ in (1.2), one finds the pole structure of
the generalized zeta function $\zeta (A,P,s)=\Tr(AP^{-s})$ (note that
the double poles become simple). In
particular, it has the Laurent expansion at zero, when $\Pi
_0(P)=0$:$$ 
\zeta (A,P,s)\sim \frac 1s {C_{-1}(A,P)}+C_0(A,P)+\sum _{l\ge
1}C_{l}(A,P)s^l,
\tag1.10
$$
where $C_0(A,P)$ is defined above and $C_{-1}(A,P)= {c'_0}$. When
$\Pi _0(P)\ne 0$, (1.10)
holds with $ C_0(A,P)$ replaced
by $ C_0(A,P)-\Tr(A\Pi _0(P))$. 

\medskip\subhead 1.2 Special cases \endsubhead

When $A=I$, we get the ordinary zeta function for $P$, $\zeta
(P,s)=\zeta (I,P,s)$. Here $c'_0=0$ and $c''_0$ is locally determined;
this also holds for $\zeta (A,P,s)$ when $A$ is a differential
operator, cf\. \cite{GS1, 
Th\. 2.7}. So $C_0(A,P)$ is locally determined when $A$ is a
differential operator. Note that when $\Pi _0(P)=0$, $\zeta
(P,0)$ is local because it equals $C_0(I,P)$. In the case $A=I$, the term
$C_1(I,P)$ in (1.10) equals minus the ``zeta-determinant''  $\log \det
P$. The higher Laurent coefficients $C_j(A,P)$ are described in
\cite{G7}.

Another interesting special case is where $A=C$ and $P=|C|$ for a
selfadjoint elliptic first-order $\psi $do $C$. Let $$
C'=C+\Pi _0(C).\tag1.11
$$
Here we can define the eta
function for $C$,$$
\eta (C,s)=\Tr(C|C|^{-s-1})=\zeta (A_C,|C|,s), \text{ with
}A_C=C|C'|^{-1}.\tag1.12 
$$  
By (1.2), it has the pole structure (note that $C$ is zero on the
nullspace of $|C|$):
$$
\eta (C,s)\sim
\frac1{\Gamma (s)}\Bigl[\sum_{ j\in \Bbb N } \frac{d_{j}}{s+{j
-n}} +  \sum_{k\in \Bbb N}\Bigl(\frac{ 
d'_{k}}{(s+k)^2}+\frac{d''_{k}}{s+ k}\Bigr)\Bigr] ,\tag1.13
$$ 
where, as usual, the division by the gamma factor makes the
double poles simple. In particular, the behavior at $s=0$ is:
$$\gathered
\eta (C,s)\sim \frac 1s
C_{-1}(A_C,|C|)+C_0(A_C,|C|)+O(s), \text{ with}\\ 
C_{-1}(A_C,|C|)= d'_0,\quad C_0(A_C,|C|)=d_n+d''_0.
\endaligned\tag1.14
$$
We recall from \cite{APS2} and Gilkey \cite{Gi} the nontrivial result that
the eta residue always vanishes:

\proclaim{Theorem 1.1} For any first-order selfadjoint elliptic $\psi
$do $C$, $d'_0=0$ in {\rm (1.13)--(1.14)}. 
\endproclaim 

Thus $\eta (C,s)$ has a finite {\it value} $d_n+d''_0$ at $0$. In
\cite{APS1}, the {\it eta invariant} was defined as this value plus
the kernel dimension of $C$;$$
\eta _C=\eta (C,0)+\dim V_0(C)=d_n+d''_0+\dim V_0(C);\tag1.15
$$
addition of other integers than $\dim V_0(C)$ can be relevant, see
(2.33)--(2.34) below.

A third special case is where $A$ is a pseudodifferential projection,
i.e., a classical $\psi $do $\Pi $ on $X$ of
order 0 satisfying $\Pi ^2=\Pi $. Here we have:

\proclaim{Theorem 1.2}

{\rm (i)} For any $\psi $do projection $\Pi $,
 $\operatorname{res}\Pi =0$; in other words, $C_{-1}(\Pi ,P)=0$ for 
arbitrary $P$. 

{\rm (ii)} If $\Pi $ is the positive eigenprojection for a
selfadjoint first-order elliptic $\psi $do $C$, then in the expansions
{\rm (1.2), (1.10)ff\.} for $A=\Pi $, $P=|C|$,$$
C_0(\Pi ,|C|)=\tfrac12 \eta (C,0)+\tfrac12 \zeta (|C|,0),\tag1.16$$
where $\zeta (|C|,0)=C_0(I,|C|)-\dim V_0(C)$; here
 $C_0(I,|C|)$ is locally determined. 

{\rm (iii)} For any orthogonal $\psi $do projection $\Pi $, there
exists a selfadjoint first-order elliptic invertible $\psi $do $C$ such that
$\Pi $ is its positive eigenprojection.
\endproclaim 

\demo{Proof} (i) was shown in \cite{W, Cor\. 7.12}. The statement
there pertains to orthogonal projections, but for any general $\psi $do
projection $\Pi $ there is an orthogonal projection $\Pi
_{\operatorname{ort}}$ with the same range and with $\Pi =R^{-1}\Pi
_{\operatorname{ort}}R$ for a suitable $\psi $do $R$ (see e.g.\
\cite{G6, Prop\. 4.8} for details and references), so
$\operatorname{res}\Pi =\operatorname{res}\Pi
_{\operatorname{ort}}=0$. 

In (ii), the positive eigenprojection means the orthogonal projection
onto the closure of the span of eigensections with positive
eigenvalues. Here $\Pi =\frac12(C+|C|)|C'|^{-1}$ (with notation
(1.11)), so in view of (1.14),$$
\aligned
\zeta (\Pi ,|C|,s)&=\tfrac12 \zeta (C|C'|^{-1},|C|,s)+\tfrac12\zeta
(I-\Pi _0(C),|C|,s)\\
&=\tfrac12 \eta (C,s)+\tfrac12\zeta (|C|,s),
\endaligned\tag1.17
$$
since $|C|^{-s}$ is taken to be 0 on $V_0(C)$.
The residue at 0 vanishes in view of Theorem 1.1 and the information
on zeta functions recalled in the beginning of this section. The
value of the left-hand 
side at zero is $\zeta (\Pi ,|C|,0)=C_0(\Pi ,|C|)$, since $\Pi \,\Pi
_0(|C|)=0$; 
and the value of the right-hand side is  $\frac12 
\eta (C,0)+\frac12\zeta (|C|,0)$.

For (iii), choose an auxiliary first-order elliptic selfadjoint $\psi
$do $C_1$ and set $C''=\Pi C_1\Pi-\Pi ^\perp C_1\Pi ^\perp$; it can be
modified on its nullspace to play the role of $C$ (more details
e.g.\ in \cite{G6, Prop.\ 4.8}).

With (iii), we see that (i) is a rather direct consequence of Theorem
1.1 and the calculations for (ii);
the argument of \cite{W} for (i) is closely related to this.\qed
\enddemo 

If in Theorem 1.2 (ii), $C$ is a differential operator, then
$C|C'|^{-1}$ is even-odd of order 0 and $C^2$ is even-even, so 
if $n$ is even, the canonical trace of $C|C'|^{-1}$ is defined
according to \cite{G7}, and
$$
\eta (C,0)=\zeta (C|C'|^{-1},|C|,0)=\zeta (C|C'|^{-1},C^2,0)=\TR
(C|C'|^{-1}).\tag1.18
$$
Then (1.16) takes the form
$$
C_0(\Pi ,|C|)=\tfrac12\TR (C|C'|^{-1})+\tfrac12\zeta (|C|,0).\tag1.19 
$$

\bigskip

\head 2. Spectral boundary
problems for first- and second-order operators\endhead 

\medskip\subhead 2.1 First-order operators \endsubhead

We now turn to boundary value problems for elliptic differential
operators, with $\psi $do
projections in the boundary condition.
The new results in this and the following chapter are published in \cite{G6}.

For first-order elliptic operators in vector bundles, it is not
always possible to get a solvable problem modulo finite dimensional
spaces (a Fredholm realization) by imposing a local boundary
condition. In their study of Dirac operators, Atiyah, Patodi and Singer \cite{APS1} therefore
considered a boundary condition with a $\psi $do projection defined from the
tangential part of the interior operator. The general choices of projections
leading to Fredholm realizations were described by Seeley in
\cite{S}; he called them {\it well-posed} boundary conditions (also
for cases of higher order operators). The well-posed problems
associated with first-order operators generalizing Dirac operators
are described in detail in \cite{G2}. 

Let us recall the set-up: 
We consider a first-order differential operator $D$
from $C^\infty (X,E)$ to $C^\infty (X,E_1)$, where $E$ and
$E_1$ are Hermitian $n_1$-dimensional vector bundles over
a compact $n$-dimensional  $C^\infty $ manifold $X$ with boundary
$\partial X=X'$. $E$ and $E_1$ have 
Hermitian metrics, and $X$ has a smooth volume element, defining 
Hilbert space structures on the sections, $L_2(E)$, $L_2(E_1)$. The
restrictions of $E$ and $E_1$ 
to the boundary $X'$ are denoted $E'$ resp\. $E'_1$. A neighborhood of
$X'$ in $X$ has the form $X_c=X'\times[0,c[\,$, and there $E,E_1$
are isomorphic to the pull-backs of $E'$ resp\. $ E'_1$. We denote
points in $X'$ resp.\ $[0,c[\,$ by $x'$ 
resp.\ $x_n$.  
$L_2(E')$  and $L_2(E'_1)$ 
are defined with respect to the volume element $v(x',0)dx'$ on $X'$
induced by the element 
$v(x',x_n)dx'dx_n$ on $X$.

$D$ may always be 
written in the following form over $X_c$:
 $$D=\sigma(\partial _{x_n}+A_1),\tag2.1
 $$
where $\sigma $ is an isomorphism from $E|_{X_c}$ to $E_1|_{X_c}$,
and $A_1$ for each $x_n$ is an elliptic operator in the $x'$-variable.
We shall study operators that resemble Dirac operators in their
structure. $D$ will be said to be of {\it product type}, when 
$\sigma $ is independent of $x_n$ and is unitary from $E'$ to
$E'_1$, and $A_1=A$ independent of $x_n$ and formally selfadjoint;
here the product measure $v(x',0)dx'dx_n$ is used on $X_c$. 
We say that $D$ is of {\it non-product type} when 
 $\sigma $ is still independent of $x_n$ and unitary, but the
condition on $A_1$ is relaxed to:
$$ A_1=A+x_nA_{11}+A_{10},\tag2.2
 $$
where $A$ is as above and
the $A_{1j}$  are smooth $x_n$-dependent tangential differential
operators in $x'$ of order $\le j$. 
In \cite{GS1}, \cite{G2}, these operators of product type and of
non-product type 
were said to be ``of Dirac-type''. In some other works, that  
notation is reserved for operators that moreover satisfy$$
E=E_1,\quad \sigma ^2=-I,\quad \sigma A=-A\sigma ,\tag2.3
$$
$D$ is formally selfadjoint on $X$ and $D^2$ is principally
scalar;
such assumptions will here only be made in special cases.

Integration by parts shows that the formal adjoint  
$D^*$ equals
$$
D^\ast=(-\partial _{x_n}+A_1')\sigma ^*,\quad A_1'=A+x_n A^*_{11}+
A'_{10}, \text{ on }X_c,\tag2.4
$$
for some morphism (zero-order operator) $A'_{10}$.

We shall also use the notation
$$
D^0=\sigma (\partial _{x_n}+A),\quad {D^0}'=(-\partial
_{x_n}+A)\sigma ^*;\tag2.5
$$these operators have a meaning on $X^0=X'\times\crp$; and ${D^0}'$ is the
formal adjoint of the operator $D^0$ going from $L_2(E^0)$ to
$L_2(E^0_1)$, where  $E^0$ and $ E^0_1$ are 
the liftings 
of  $E'$ resp\. $E'_1$ to $X^0$, and the product measure is used.

Sometimes we consider, along with $D$ satisfying (2.1)--(2.2), a product
type operator $D_0$ on $X$, which is
like $D^0$ on $X_c$ and is extended to 
an elliptic operator on $X$ (e.g\. by patching it together with $D$).

By $V_>$, $V_\geq$, $V_<$ or $V_\leq$  we denote the subspaces of $L^2(E')$
spanned by the eigenvectors of $A$ corresponding to eigenvalues which
are $>0$, $\geq0$, $<0$, or $\leq0$. (For precision one can write
$V_>(A)$, etc.) $V_0$
is the nullspace of $A$. 
The corresponding projections are denoted $\Pi_>$, $\Pi _\ge$, etc\. (note
that $\Pi _\ge=\Pi _>+\Pi _0$ and $\Pi _<=I-\Pi _\ge$). They are
pseudodifferential operators 
($\psi $do's) of order 0; $\Pi _0$ has finite rank and is a $\psi $do
of order $-\infty $. 
We set
$$\gathered
|A|=(A^2)^\frac12,\quad A'=A+\Pi _0,\quad \text{ so that }|A'|=|A|+\Pi
_0\text{ and}\\
\Pi_>=\tfrac12\tfrac{|A|+ A}{|A'|}=\tfrac12\tfrac A
{|A'|}+\tfrac12-\tfrac12 \Pi _0.
\endgathered\tag2.7
$$

Together with the equation $Du=f$, we consider a
boundary condition$$
\Pi \gamma _0u=0,\tag2.8
$$
 where $\gamma_0u=u|_{X'}$, defining the realizations $D_{\Pi }$ and
$D^0_\Pi $, acting like $D$ resp\. $D^0$ and with domain $$
D(D_\Pi )\text{ resp\. }D(D^0_{\Pi })=\{\, u\in H^1(E)\text{ resp\.
}H^1(E^0)\mid \Pi \gamma _0 u=0\,\};\tag2.9 
$$
we denote by $H^s(E)$ the Sobolev space of order $s$. \cite{APS1}
considered the case where $\Pi =\Pi _\ge(A)$, 
but increasingly general projections have been studied through
the years. The most general case is where
$\Pi $ is a pseudodifferential projection that is {\it
well-posed} with respect to $D$ (cf\. Seeley \cite{S} or \cite{G2}).
This means 
that when we at each $(x',\xi ')$ in the cotangent sphere bundle of
$X'$ denote by $N^+(x',\xi ')\subset \Bbb C^{n_1}$ the space
of boundary values of null-solutions of the model operator (defined
from the principal symbol $d^0$ of $D$ at $X'$),$$
N^+(x',\xi ')=\{\, z(0)\in \Bbb C^{n_1}\mid d^0(x',0,\xi ',D_{x_n})z(x_n)=0
,\; z(x_n)\in L_2(\Bbb R_+)^{n_1}\,\},\tag2.10
$$
then the
principal symbol $\pi ^0(x',\xi ')$ of $\Pi $ maps $N^+(x',\xi ')$
bijectively onto the range of $\pi ^0(x',\xi ')$ in $\Bbb C^{n_1}$. 
Equivalently, the model problem with homogeneous boundary condition
is uniquely solvable in $L_2(\Bbb R_+)^{n_1}$.

\example{Example 2.1}
For $|\xi '|\ge 1$, the space $N^+(x',\xi ')$ equals the positive
eigenspace for $a^0(x',\xi ')$, i.e., the range of the principal
symbol $\pi ^0_>(x',\xi ')$ of $\Pi _\ge(A)$, so $\Pi _\ge (A)$
is well-posed for $D$; this is the case considered in \cite{APS1}.
Various finite rank perturbations of $\Pi _\ge (A)$ were considered in
Douglas-Woj\-cie\-chow\-ski \cite{DW}, M\"uller \cite{M}, Dai and
Freed \cite{DF}, Grubb and Seeley \cite{GS1, GS2}. Booss-Bavnbek and
Wojciechowski, cf\. e.g. \cite{BW}, pointed to the interest of
studying the exact Calder\'on projector 
which differs from $\Pi _\ge (A)$ by an operator of order $-\infty $
in the product case; \cite{Woj} treated quite general
perturbations of order $-\infty $. Br\"uning and Lesch
\cite{BL} studied a principally different family of pseudodifferential
projections that we shall here denote $\Pi (\theta )$, and finally
\cite{G2, G4} included all well-posed projections in the study.
\endexample

\medskip\subhead 2.2 Second-order operators  \endsubhead

As noted in \cite{S}, \cite {G2} (and in the proof of Theorem 1.2
(i) above), it is no restriction to assume that $\Pi $ is an
orthogonal projection. In view of Green's formula
$$
(Du,v)_{X}-(u,D^*v)_{X}=
-(\sigma \gamma _0u,\gamma _0v)_{X'},\tag2.11
$$
and elliptic regularity, the adjoint $(D_\Pi )^*$ is the realization of $D^*$
defined by the
boundary condition $\Pi ^\perp\sigma ^*\gamma _0v=0$ (associated with
the well-posed projection $\Pi '=\sigma \Pi ^\perp\sigma ^*$ for $D^*$).
It follows that $D^*D$ is of the form (on $X_c$):
$$
P=-\partial _{x_n}^2+P'+x_nP_2+P_1,\tag2.12
$$ 
with $P'=A^2$, the $P_j$ being $x_n$-dependent
differential operators of order $j$ in $E|_{X_c}$, and that
${D_\Pi }^*D_{\Pi }$ 
is the realization of $D^*D$ defined by the
boundary condition $$
\Pi \gamma _0u=0,\quad \Pi ^\perp(\gamma _1u+A_1(0)\gamma _0u)=0.\tag2.13
$$ 
The study of spectral invariants of $D_\Pi $ can to a large extent
be based on the study of the second-order realization ${D_\Pi
}^*D_{\Pi }$.

Recently, there has been an interest in studying similar second-order
problems for their own sake, with a view to applications in brane
theory (see Vassiliev \cite{V1}, \cite{V2} and the introduction in
\cite {G6}), so let us look at a slightly more general
situation:

$P$ is an elliptic second-order partial differential operator in
$E$, of the form  (2.12) on $X_c$, with $P'$ being an elliptic selfadjoint
nonnegative second-order
differential operator in $E'$ (independent of $x_n$), and the $P_j$
as described above. It is considered together with the boundary condition
$$
Tu=0,\quad\text{where }Tu=\{\Pi _1\gamma _0u, \Pi _2(\gamma
_1u+B\gamma _0u)\},\tag2.14
$$
where $\Pi _1$ is a $\psi $do projection operator and $B$ is
a first-order $\psi $do, both acting in $E'$, and $\Pi _2=I-\Pi _1$.
We denote by $P_T$ the realization of $P$ defined by this boundary
condition; it acts like $P$ and has the domain
$$
D(P_T )=\{\, u\in H^2(E)\mid T u=0\,\}.\tag2.15
$$

Here ${D_\Pi }^*D_{\Pi }=P_T$ in the special case where $P=D^*D$,
$\Pi _1=\Pi $ and
is orthogonal, and $B=A_1(0)$.

We denote $\frak A=(P'-\lambda )^{\frac12}=(P'+\mu ^2)^{\frac12}$; here $\lambda $ runs in $\Bbb C\setminus\crp$, and
$\mu =(-\lambda )^{\frac12}$ runs in $\{ \mu \mid
\operatorname{Re}\mu >0\}$. The principal symbol is $\frak a^0(x',\xi
',\mu )=(p^{\prime 0}(x',\xi ')+\mu ^2)^\frac12$.
We can assume that $X$ is smoothly imbedded in a closed
$n$-dimensional manifold $\widetilde X$, provided with a vector
bundle $\widetilde E$ such that $E=\widetilde E|_X$,
and such that $P$ is
defined in $\widetilde E$ with similar properties.  
Let us denote 
$$
\Gamma _{\theta }=\{\, \mu \in\Bbb C\setminus \{0\}\mid |\arg\mu
|<\theta\,\}.\tag2.16
$$
The following result was proved in \cite{G6, Sect\. 2}:

\proclaim{Theorem 2.2} Assume {\rm (H1)} and {\rm (H2)}:
\roster
\item"{\rm (H1)}" The principal symbols of $\Pi _1$ and $P'$ commute.
\item"{\rm (H2)}" There is a $\theta \in \,]0,\frac\pi 2]$ such
that, with  $b ^h(x',\xi  ')$ and  $\pi _i ^h (x',\xi 
')$  denoting the strictly homogeneous principal
symbols of 
$B$ and the $\Pi _i$,$$
\frak a^0 -\pi _2^h b^h  \pi ^h
_2\text{ is
invertible  for }\xi '\in\Bbb R^n,\mu \in\Gamma _\theta\cup \{0\}\text{ with }(\xi ',\mu
)\ne (0,0).\tag2.17
$$ 
\endroster
Then for
each $\theta '\in \,]0,\theta [\,$ there is an $r=r(\theta ')\ge 0$ such
that when $|\arg\mu |\le \theta '$ and $|\mu |\ge  r$,
$P_T+\mu ^2=P_T-\lambda $ is a bijection from 
$D(P_T)$ to $L_2(E)$ with inverse $(P_T-\lambda )^{-1}=R_T(\lambda )$;
$$
R_T(\lambda )=Q(\lambda )_++G(\lambda ),\tag2.18
$$
where $Q(\lambda )=(P-\lambda )^{-1}$ on $\widetilde X$ and $G(\lambda )$
is a singular Green operator belonging to the parameter-dependent
calculus of \cite{G3}, 
with symbol in
$\Cal
S^{0,0,-3}(\Gamma _{\theta },\Cal 
S_{++})$.
\endproclaim 
\comment
$$
\aligned
R_T(\lambda )&=Q(\lambda )_++G(\lambda ),\\
G(\lambda )&=-K_{\operatorname{D}}(\lambda )\gamma _0Q(\lambda )_++K_{\operatorname{D}}(\lambda ) [S_0(\lambda )\gamma
_0+S_1(\lambda )\gamma _1]Q(\lambda )_+;
\endaligned\tag2.35
$$ 
here $S_0$ and $S_1$ (given in {\rm (2.36)} below) are weakly
polyhomogeneous $\psi $do's in $E'$ 
lying in $\operatorname{OP}'S^{0,0,0}(\Gamma _{\theta })$ 
resp\. $\operatorname{OP}'S^{0,0,-1}(\Gamma _{\theta })$,  and hence
$G$ is a singular Green operator of class $0$ 
\endcomment

The original result also gave a more precise formula for $G(\lambda )$ in
terms of the given operators; we shall consider some
consequences of this later.

Using the general machinery of \cite{G3} (or more specific calculations as in \cite{GS1}, \cite{G2}), one
deduces the existence of the following general trace expansions:

\proclaim{Theorem 2.3} Assumptions as in Theorem {\rm 2.2}.
Let $F$ be a differential operator in $E$ of order
$m$ and let $N>\frac{n+m}2$. Then $FR^N_T(\lambda )$ is trace-class and
the trace has an expansion for $|\lambda |\to\infty $ with $\arg
\lambda \in \,]\pi -2\theta , \pi +2\theta [\,$ (uniformly in closed
subsectors):$$
\multline\Tr\bigl( F R_T^N(\lambda )\bigr)\\
\sim
\sum_{-n\le k<0 } \tilde a_{ k}(F)(-\lambda ) ^{\frac{{m} -k}2-N}+ 
\sum_{k\ge 0}\bigl({ \tilde a'_{ k}(F)}\log (-\lambda ) +{\tilde a''_{
k}(F)}\bigr)(-\lambda ) ^{\frac{{m} -k}2-N},
\endmultline\tag2.19
$$
with locally determined coefficients $\tilde a_k$ and $\tilde a'_k$.
If
$m$ is odd, 
$\tilde a_{-n}=0$.

Here, if $F$ is tangential (differentiates only with respect to $x'$)
on $X_c$, the log-coefficients $\tilde a'_k$ with $0\le k<m$ vanish,
and the $\tilde 
a''_k$ with $0\le k<m$ are locally determined.
\endproclaim 

In these formulas, the notation differs slightly from that of (1.1):
we have collected the coefficients of each power
of $-\lambda $ in one term, and adapted the indexation of all
coefficients to the way $k$ enters in the powers. 

The expansion can be translated (as in \cite{GS2}) to a statement on
the meromorphic 
extension of $\Tr(FP_T^{-s})$; let us write the result in the cases
where $F=\varphi $, a morphism (or ``smearing function''), or $F=D_1$, a
first-order differential operator:
$$\aligned
&\Gamma (s)\Tr (\varphi P_T^{-s})\\
&\sim
\sum_{-n\le k<0}\frac {a_{k}(\varphi )}{s+\frac k2}
-\frac{\Tr(\varphi  \Pi _0(P_T))}s+\sum_{k=0}^\infty
\Bigl(\frac {a'_{k}(\varphi )}{(s+\frac k2)^2}+\frac{a''_{k}(\varphi
)}{s+\frac k2}\Bigr);\\
&\Gamma (s)\Tr (D_1 P_T^{-s})\\
&\sim
\sum_{-n< k<0}\frac {a_{k}(D_1)}{s+ \frac{k-1}2}
-\frac{\Tr(D_1  \Pi _0(P_T))}{s}+\sum_{k=0}^\infty
\Bigl(\frac {a'_{k}(D_1)}{(s+ \frac{k-1}2)^2}+\frac{a''_{k}(D_1)}{s+ \frac{k-1}2}\Bigr);
\endaligned\tag2.20$$
the last expansion can also be written in the more customary form
(with $s=\frac{s'+1}2$):
$$
\multline
\Tr (D_1 P_T^{-\frac{s'+1}2})\\
\sim
\frac 1{\Gamma (\frac{s'+1}2)}\Bigl[\sum_{-n< k<0}\frac {2a_{k}(D_1)}{s'+ k}
-\frac{2\Tr(D_1  \Pi _0(P_T))}{s'+1}+\sum_{k=0}^\infty
\Bigl(\frac {4a'_{k}(D_1)}{(s'+ k)^2}+\frac{2a''_{k}(D_1)}{s'+
k}\Bigr)\Bigr].
\endmultline \tag2.21
$$
(The formula (1.13) was simpler, since $|C'|^{-1}$ could be taken
into the operator in front.)
The coefficients are related to those in the resolvent expansions by
universal nonzero proportionality factors; in particular,
$$
\tilde a'_0(F)= a'_0(F);\quad \tilde a''_0(F)=a''_0(F).\tag2.22
$$

If $\theta >\frac\pi 4$ in (H2) so that the ``heat
operator'' $e^{-tP_T}$ exists, there is a trace expansion
of $Fe^{-tP_T}$ in the spirit of (1.3).

In (2.20),  $\Tr (\varphi P_T^{-s})$ is a generalized zeta function,
also denoted $\zeta (\varphi ,P_T,s)$, and $\Tr (D_1 P_T^{-\frac{s+1}2})$
in (2.21) is somewhat like an eta function. 

\medskip\subhead 2.3 Consequences for first-order operators \endsubhead

It is accounted for in \cite{G6} how 
$P_T={D_\Pi }^*D_{\Pi }$ enters as a special case in the above theorems,
when the principal 
symbols of $A^2$ and $\Pi $ commute (this holds in particular when $A^2$ is
principally scalar, i.e., the principal symbol is scalar). Here the
well-posedness of $\Pi $ is in a certain sense equivalent with (H2).
So, our basic assumptions on $D$ and $\Pi $ are as follows:

\proclaim{Hypothesis (H3)} $D$ is as described in {\rm (2.1)}ff., of
product type or non-product type. $\Pi $ is an orthogonal $\psi $do
projection in $L_2(E')$ that is well-posed for $D$, and the principal
symbols of $\Pi $ and $A^2$ commute.
\endproclaim 

Then $\zeta (\varphi ,P_T,s)$ equals $ \zeta (\varphi ,{D_\Pi }^*D_{\Pi },s)$,
also denoted $\zeta ({D_\Pi }^*D_{\Pi },s)$ if $\varphi =I$.
We can let
$D_1=\psi D$ for some morphism $\psi $
from $E_1$ to $E$, defining 
$$
\Tr (D_1 P_T^{-\frac{s+1}2})=\Tr(\psi D ({D_\Pi }^*D_{\Pi
})^{-\frac{s+1}2})=\eta (\psi ,D_\Pi ,s),\tag2.23
$$
an eta function of $D_\Pi $.
The existence of the above expansions for \linebreak$\Tr(\varphi ({D_\Pi
}^*D_{\Pi }-\lambda )^{-N})$, $\Tr(\psi D({D_\Pi
}^*D_{\Pi }-\lambda )^{-N})$, 
$\zeta (\varphi ,{D_\Pi }^*D_{\Pi
},s)$ and $\eta 
(\psi ,D_\Pi ,s)$ is known from
\cite{G2} for general $\Pi $ (with special choices of $\Pi $
treated in earlier works); we repeat them here for clarity:$$\aligned
&\Tr( \varphi (D_\Pi^*D_{\Pi }-\lambda )^{-N})\\
&\quad\sim
\sum_{-n\le k<0 } \tilde a_{ k}(\varphi )(-\lambda ) ^{-\frac{ k}2-N}+ 
\sum_{k\ge 0}\bigl({ \tilde a'_{ k }(\varphi )}\log (-\lambda ) +{\tilde a''_{
k}(\varphi )}\bigr)(-\lambda ) ^{-\frac{k}2-N},\\
&\Tr(D_1(D_\Pi
^*D_{\Pi }-\lambda )^{-N}) \\
&\quad\sim\sum_{-n< k<0 } \tilde a_{ k}(D_1)(-\lambda ) ^{\frac{1 -k}2-N}+ 
\sum_{k\ge 0}\bigl({ \tilde a'_{ k}(D_1)}\log (-\lambda ) +{\tilde a''_{
k}(D_1)}\bigr)(-\lambda ) ^{\frac{1 -k}2-N},\\
\endaligned$$
$$\aligned
&\Gamma (s)\Tr (\varphi (D^*_\Pi D_\Pi )^{-s})\\
&\quad\sim
\sum_{-n\le k<0}\frac {a_{k}(\varphi )}{s+\frac k2}
-\frac{\Tr(\varphi  \Pi _0(D_{\Pi }))}s+\sum_{k=0}^\infty
\Bigl(\frac {a'_{k}(\varphi )}{(s+\frac k2)^2}+\frac{a''_{k}(\varphi
)}{s+\frac k2}\Bigr),\\
&\Tr (D_1 (D^*_\Pi D_\Pi )^{-\frac{s+1}2})\\
&\quad\sim
\frac 1{\Gamma (\frac{s+1}2)}\Bigl[\sum_{-n< k<0}\frac {2a_{k}(D_1)}{s+ k}
+\sum_{k=0}^\infty
\Bigl(\frac {4a'_{k}(D_1)}{(s+ k)^2}+\frac{2a''_{k}(D_1)}{s+
k}\Bigr)\Bigr];
\endaligned\tag2.24$$
where we used that $D_1=\psi D$ vanishes on $V_0(D_\Pi )$.

The coefficient analysis for
$P_T$ described 
below will allow some new conclusions on these special cases also.

\medskip\subhead 2.4 Analysis of the zero'th coefficients \endsubhead

The formulas (2.20) and (2.21) show in particular how the zeta function and
eta-like function 
behave near $s=0$: $$\aligned
\Tr (\varphi P_T^{-s})&= a'_0(\varphi )s^{-1}+(a''_0(\varphi )-\Tr(\varphi \Pi
_0(P_T))s^0+O(s),\\
\Tr (D_1 P_T^{-\frac{s+1}2})&=\pi ^{-\frac12} 2a'_0(D_1)s^{-2}+\pi ^{-\frac12} 4a''_0(D_1)s^{-1}
+O(1),
\endaligned\tag2.25$$
for $s\to 0$.
An important question in this context is what we can say about the
value, or the vanishing, of the coefficients in the Laurent
expansions (2.25). It is here that we use the more precise
description of the singular Green part $G(\lambda )$ of $R_T(\lambda
)$ mentioned after Theorem 2.2.

We recall from the general theory of pseudodifferential boundary
operators ($\psi $d\-bo's) that when a singular Green operator
$G=\operatorname{OPG}(g(x',\xi ,\eta _n,\mu ))$ is trace-class on
$\rnp$, its trace equals the $\Bbb R^{n-1}$-trace of the $\psi $do on
$\Bbb R^{n-1}$ called the {\it normal trace} of $G$, $\tr_nG$; it is the
operator with symbol
$$
(\operatorname{tr}_n g)(x',\xi ',\mu )=\int g(x',\xi ',\xi
_n,\xi _n,\mu
)\,\d\xi _n.$$
Then the trace expansion of $G$ is obtained by applying the rules for
the boundaryless manifold $\Bbb R^{n-1}$ to $\tr_nG$.
We also observe that
$$
R^N_T=(Q^N)_++G^{(N)}=\tfrac1{(N-1)!}{\partial _\lambda ^{N-1}}
R_T
=\tfrac1{(N-1)!}{\partial _\lambda ^{N-1}}Q_+ +
\tfrac1{(N-1)!}{\partial _\lambda ^{N-1}}G,
\tag2.26$$
where 
$G^{(N)}=\frac1{(N-1)!}\partial _\lambda ^{N-1}G$ is a singular
Green operator of class $0$  
with symbol in $\Cal S^{0,0,-2N-1}(\Gamma _{\theta },\Cal
S_{++})$.

The $\psi $do $\varphi Q^N(\lambda )_+$ has a trace expansion without
logarithmic or nonlocal terms:$$
\Tr \varphi Q^N(\lambda )_+\sim\sum_{k\ge -n, k+n\text{
even}}c_{k}(\varphi )(-\lambda
)^{-\frac k2 -N}.\tag 2.27
$$

The crucial information on $G(\lambda )$ that we shall use is shown in
\cite{G6, Sect\. 4}:

\proclaim{Theorem 2.4} 
 Let $\varphi $ be a morphism in $E$, independent of $x_n$ on $X_c$ (its restriction to $X'$ likewise
denoted $\varphi $). For $G^{(N)}(\lambda )$ from {\rm (2.26)},
cut down to $X_c$, we have that$$
\tr_n\varphi G^{(N)}(\lambda )= \tfrac12\varphi \Pi _2 (P'-\lambda
)^{-N}+ S_1(\lambda ) +S_2 (\lambda ),\tag2.28
$$
where $S_1$ and $S_2$ have trace expansions of the form$$\aligned
\Tr_{X'}S_1(\lambda )&\sim\sum_{k\ge 1 -n}s_{1,k}(-\lambda )^{-\frac k2 -N},\\
\Tr_{X'}S_2(\lambda )&\sim\sum_{1-n\le k\le 0}s_{2,k}(-\lambda
)^{-\frac k2 -N} 
+\sum_{k\ge 1}(s'_{2,k}\log(-\lambda )+s''_{2,k})(-\lambda )^{-\frac k2 -N};
\endaligned\tag2.29$$
here the $s_{i,k}$ and $s'_{i,k}$ are locally determined.
\endproclaim  

On interior coordinate patches, the trace of $G^{(N)}(\lambda )$ is
$O(\lambda ^{-M})$, for any $M$. Thus the only contributions
to $a'_0(\varphi )$ and the only nonlocal contributions to
$a''_0(\varphi )$ in (2.20), (2.25) come from the 
first term in 
the right-hand side of (2.28)! And this is a function whose expansion
we know very well from the case of closed manifolds.
In fact, $\tfrac12\varphi \Pi _2 (P'-\lambda
)^{-N}$ has an expansion:
$$\multline
\Tr\bigl(\tfrac12\varphi \Pi _2 (P'-\lambda )^{-N}\bigr)
\\ \sim
\sum_{ 1-n\le k<0  } \tilde c_{ k}(-\lambda ) ^{\frac k2-N}+ 
\sum_{k\ge 0, k\text{ even}}\bigl( \tilde c'_{ k}\log (-\lambda )
+\tilde c''_{ k} \bigr)(-\lambda ) ^{-\frac k2-N},\endmultline
\tag2.30$$
by (1.1) (with $A=\tfrac12\varphi \Pi _2 $ of order 0,  $\dim X'=n-1$,
and a regrouping and change in the indexation as indicated after (2.19)).
Here, since $\operatorname{res}\varphi =0$, $$
\tilde c'_0=\tfrac12\operatorname{res}(\tfrac12\varphi \Pi _2
)=-\tfrac14\operatorname{res}(\varphi \Pi _1 ).\tag2.31
$$
We conclude immediately, in view of Theorem 1.2 (i):

\proclaim{Theorem 2.5} Assumptions of Theorem {\rm 2.3}. One has in
general that
$a'_0(I)=0$, and $a'_0(\varphi )=-\frac14\res(\varphi \Pi _1)$, in
{\rm (2.25)}.  

In particular, $\zeta ({D_\Pi }^*D_{\Pi },s)$ is
regular at $s=0$ for all choices of $\Pi $.
\endproclaim 

Note that $B$ does not enter in the value. Also, the result for
$\zeta ({D_\Pi }^*D_{\Pi },s)$ shows that the regularity 
at $s=0$ is preserved under perturbations of $\Pi $ by
operators of order $\le -1$, since such perturbations preserve
well-posedness. (This was known earlier for perturbations of order
$\le -n$, \cite{G2}.)

\example{Remark 2.6}
When $\varphi $ is nontrivial, there is another sufficient condition for
the vanishing of $a'_0(\varphi )$ (apart from the possibility that $\varphi
\Pi _1$ could be a 
projection): When $\Pi _1$ is the positive 
eigenprojection for a 
selfadjoint {\it differential operator} $C$ and $n-1$ is even, then
$$
\operatorname{res}(\varphi \Pi _1)=\operatorname{res}(\tfrac12
\varphi (I+C|C'|^{-1}))=\tfrac12\operatorname{res}(\varphi C|C'|^{-1})
$$
vanishes since $\varphi C|C'|^{-1}$ has even-odd parity, cf.\ (1.8).
Also perturbations of $\Pi _1$ of order
$\le -n$ are allowed, since they do not interfere with the residue.
\endexample 

\example{Remark 2.7} Concerning $\tilde c''_0$, we observe: When $\Pi
_1$ is the positive eigenprojection 
for a first-order selfadjoint invertible elliptic $\psi $do $C$, then
in (2.30) with $\varphi =I$,$$
\aligned
\tilde c''_0&=C_0(\tfrac12\Pi _2, P')=-\tfrac12 C_0(\Pi _1, P')+\text{
local terms }\\
&=-\tfrac12 C_0(\Pi _1, |C|)+\text{
local terms }
=-\tfrac14 \eta (C,0)+\text{ local terms,}
\endaligned \tag2.32
$$
by (1.9), (1.16) and the fact that $\zeta (|C|,0)=C_0(I, |C|)
$ is local. The case where $C$ has a nontrivial nullspace $V_0(C)$
 is analyzed in
\cite {G6, Th\. 4.9, Cor.\ 5.4--5.5}. Here it is found e.g.\ that if
$V_0(C)=V'_0\oplus V''_0$ (orthogonal decomposition) and $\Pi _1=\Pi
_>(C)+\Pi _{V'_0}$ (the latter denoting the orthogonal projection
onto $V'_0$), then $$\aligned
\zeta (P_T,0)&=-\tfrac14 \eta _{C, V'_0}-\dim
V_0(P_T)+\text{ local contributions,}\\
\text{with }\eta _{C,V'_0}&=\eta (C,0)+\dim V'_0-\dim V''_0.
\endaligned\tag 2.33
$$
Moreover, one has for $D_\Pi $ with such $\Pi $, that$$
\operatorname{index} D_\Pi =-\tfrac12 \eta _{C, V'_0}+\text{ local
contributions,}\tag 2.34
$$
which allows the remarkable observation that the ``non-locality''
depends only on the 
projection, not on the interior operator. 
\endexample

Similarly to Theorem 2.4, one has for $D_1G(\lambda )$: 

\proclaim{Theorem 2.8} Let  
$D_1 $ be a first-order differential operator on $X$, of the form
$D_1=\psi (\partial _{x_n}+B_1)$ on $X_c$, where $B_1$ is tangential
and $\psi $ is a morphism in $E$, independent of $x_n$. Then for
$G^{(N)}(\lambda )$,
cut down to $X_c$, we have that$$
\tr_n(D_1 G^{(N)}(\lambda ))= -\tfrac12\psi  \Pi _2
\tfrac1{(N-1)!}{\partial _\lambda ^{N-1}}(P'-\lambda )^{-\frac12}+
\widetilde S_1(\lambda )
+\widetilde S_2 (\lambda ),\tag2.35
$$
where $\widetilde S_1$ and $\widetilde S_2$ have trace expansions of
the form:$$\aligned 
\Tr_{X'}\widetilde S_1(\lambda )&\sim\sum_{k\ge 1-n}\tilde
s_{1,k}(-\lambda )^{\frac {1-k}2 -N},\\ 
\Tr_{X'}\widetilde S_2(\lambda )&\sim\sum_{1-n\le k\le 0}\tilde
s_{2,k}(-\lambda )^{\frac {1-k}2 -N} 
+\sum_{k\ge 1}(\tilde s'_{2,k}\log(-\lambda )+\tilde
s''_{2,k})(-\lambda )^{\frac {1-k}2 -N}; 
\endaligned\tag2.36$$
here the $\tilde s_{i,k}$ and $\tilde s'_{i,k}$ are locally determined.
\endproclaim  

An application of the calculus of \cite{GS1} gives that
$$
\multline
\Tr_{X'}(-\tfrac12\psi  \Pi _2
\tfrac1{(N-1)!}{\partial _\lambda ^{N-1}}(P'-\lambda )^{-\frac12})\\
\sim\sum_{1-n\le k< 0}\tilde d_k(\psi )(-\lambda )^{\frac {1-k}2 -N}
+\sum_{k\ge 0}(\tilde d'_{k}\log(-\lambda )+\tilde d''_{k})(-\lambda
)^{\frac {1-k}2 -N},\endmultline\tag2.37
$$
where an analysis as in \cite{G5, pf\. of Th\. 5.2} shows that $$
\tilde d'_0=-\alpha  \operatorname{res}(\psi \Pi _2)=\alpha  \operatorname{res}(\psi \Pi _1),\tag2.38
$$
with a universal nonzero factor $\alpha $.
This is the only contribution to $a'_0(D_1)$, so we conclude:

\proclaim{Theorem 2.9} Assumptions as in Theorem {\rm 2.7}.
In {\rm (2.25)}, $$
a'_0(D_1)=\alpha \res (\psi \Pi _1),\tag2.39
$$
with a universal nonzero factor $\alpha $. Here $a'_0(D_1)$ vanishes if $\psi
\Pi _1$ is a projection. 

In particular, $\eta (\psi ,D_{\Pi },s)$ has a
simple pole at $0$ if $\psi  \sigma \Pi $ is a projection, e.g.\ if
$\psi =\sigma ^*$.  
\endproclaim 

As in Remark 2.6, another sufficient condition for the
vanishing of $a'_0(D_1)$ is that $\Pi _1$ is the positive
eigenprojection of a differential operator and $n$ is odd. Also
perturbations of $\Pi _1$ by operators of order $\le -n$ are allowed here.

\bigskip

\head 3. Results under further symmetry conditions,
perturbations of the boundary
projection \endhead 

\medskip\subhead 3.1 Results for zeta functions \endsubhead

In the following, we take $\Pi _1$ equal to an orthogonal
pseudodifferential projection $\Pi $ (so that $\Pi _2=\Pi ^\perp$).
We consider the case where there exists a unitary 
morphism
$\sigma $ in $E$ such that $$
\sigma ^2=-I,\quad \sigma P'=P'\sigma ,\quad \Pi ^\perp=-\sigma \Pi
\sigma .\tag3.1
$$

\proclaim{Theorem 3.1} Let {\rm (3.1)} hold.
Then$$\aligned
\Tr_{X'}(\tfrac12\Pi ^\perp\tfrac{\partial _\lambda
^{m-1}}{(m-1)!}(P'-\lambda )^{-1})&= \tfrac14\Tr_{X'}(\tfrac{\partial _\lambda
^{m-1}}{(m-1)!}(P'-\lambda )^{-1}),\\
-\Tr_{X'}(\tfrac12\sigma \Pi ^\perp\tfrac{\partial _\lambda
^{m-1}}{(m-1)!}(P'-\lambda )^{-\frac12})&= -\tfrac14\Tr_{X'}(\sigma
\tfrac{\partial _\lambda 
^{m-1}}{(m-1)!}(P'-\lambda )^{-\frac12}).
\endaligned\tag3.2$$
Let $\Pi _1=\Pi $,
$\Pi _2=\Pi ^\perp$.
Then in {\rm (2.30)} with $\varphi =I$, and in {\rm (2.37)} with
$\psi =\sigma $, all log-terms vanish and all the remaining
coefficients are locally determined. In particular,
$\tilde c''_0$ is locally determined (from the
symbol of $P'$), $\tilde d'_0$ vanishes, and $
 \tilde d''_0$ is locally determined
(from the symbol of $P'$ and $\sigma $).

It follows that in {\rm (2.20), (2.25)} with $\varphi =I$ and
$D_1=\sigma (\partial _{x_n}+B_1)$, $$
\gathered \tilde a'_0(D_1)= a'_0(D_1)=0,\\
\tilde a''_0(I ),a''_0(I),\tilde a''_0(D_1)\text{ and } a''_0(D_1)
\text{ are
locally determined}.
\endgathered
\tag3.3
$$
\endproclaim 

In the proof, the identities in (3.2) are obtained by linearity and
cyclic permutation in the trace formulas. Now since the operators
$(P'-\lambda )^{-a}=(P'+\mu ^2)^{-a}$, $a\in \Bbb N$ or $\Bbb
N+\frac12$, are strongly polyhomogeneous in $(\xi ',\mu )$, the
traces have expansions 
without logs and with only local coefficients, by \cite{GS1}.

More precisely, $a''_0(I)$ in this case depends on the symbol of $P$,
on $\sigma $,
and on the first $n$ strictly homogeneous terms in the symbols of $\Pi $ and
$B$; and $a''_0(D_1)$ depends on the mentioned symbols together
with that of $B_1$.

We shall pursue this result for the traces arising from $D_{\Pi }$ in
cases with selfadjointness properties. 
Here we are interested in truly selfadjoint product cases as well as
in nonproduct cases where $D$ is principally 
selfadjoint at $X'$. Assume that $E=E_1$.
Along with $D$ we consider a product type operator $D_0$, defined as
after (2.5).

In addition to
the requirements that $\sigma $ be unitary and $A$ be selfadjoint, we
now assume (2.3), which means that $D_0$ is formally selfadjoint on
$X_c$ when this is provided with the product volume element
$v(x',0)dx'dx_n$.  (If $D_0$
is selfadjoint on $X$, we call this a {\it selfadjoint product case}.)

When $\Pi $ is an orthogonal projection in $L_2(E')$, it
is well-posed for $D$ if and only if it is so for $D_0$.
For $D_0$ in selfadjoint product cases,
some choices of $\Pi $ will lead to selfadjoint realizations
$D_{0,\Pi }$, namely 
(in view of (2.11)) those for which
$$
\Pi =-\sigma \Pi ^\perp \sigma .\tag 3.4
$$
The properties (2.3) and (3.4) imply (3.1) with $P'=A^2$, so we can
apply Theorem 
3.1 to  $D_\Pi ^*D_{\Pi }$ (and $D_{0,\Pi }^2$).

As pointed out in the appendix A.1 of Douglas and Wojciechowski
\cite{DW}, it follows from Ch\. 17 (by Palais and Seeley) of the
Palais seminar \cite{P}
that when (2.3) holds and $n$ is odd, there exists a subspace $L$
of $V_0(A)$ such that $\sigma L\perp L$ and $V_0(A)= L\oplus\sigma
L$. M\"uller showed in \cite{M} (cf\. (1.6)ff\. and Prop\.
4.26 there) that such $L$ can be found in any dimension. Denoting the
orthogonal projection onto $L$ by $\Pi _L$, we have that  
$$
\Pi _+=\Pi _>(A)+\Pi _L\tag3.5
$$
satisfies (3.4). The projections $\Pi (\theta )$
introduced by Br\"uning and Lesch \cite{BL} likewise satisfy (3.4).
These projections commute with $A$, so Hypothesis (H3) is satisfied.

Theorem 3.1 implies immediately:

\proclaim{Corollary 3.2} When $D$ and $\Pi $ satisfy {\rm (H3)} and
in addition {\rm (2.3)} and {\rm (3.4)}, then in {\rm
(2.24)} with $\varphi =I$,
$$
\tilde a''_0(I) \;(=a''_0(I)) \text{ is
locally determined.}\tag3.6
$$
\endproclaim 
 
This has an interesting consequence for perturbations of $\Pi $:

\proclaim{Theorem 3.3} In addition to the hypotheses of Corollary
{\rm 3.2}, assume that$$
\Pi =\overline\Pi +\Cal S,\tag3.7
$$
where $\overline\Pi $ is a fixed well-posed projection satisfying
{\rm (3.4)} and $\Cal S$ is of order $\le 
-n$. ($\overline \Pi $ can in particular be taken as $\Pi _+$ in {\rm
(3.5)} or $\Pi (\theta )$ from \cite{BL}.)

Then the $\tilde a''_0(I)$-terms (and $a''_0(I)$-terms) in {\rm
(2.24)} for $D^*_{\overline\Pi }D_{\overline\Pi } 
$ and 
$D^*_\Pi D_\Pi 
$ are the same,$$
\tilde a''_0(I)(D_\Pi ^*D_\Pi )=\tilde a''_0(I)
(D_{\overline\Pi }^*D_{\overline\Pi }).\tag3.8
$$
It follows that
$$
\zeta (D_\Pi ^*D_\Pi  ,0)+\operatorname{dim}V_0(D_\Pi )=\zeta 
(D_{\overline\Pi }^*D_{\overline\Pi },0)+\operatorname{dim}V_0(D_{\overline\Pi } );\tag3.9
$$
in particular$$
\zeta (D_\Pi ^*D_\Pi ,0)=\zeta (D_{\overline\Pi }^*D_{\overline\Pi },0) \quad (\operatorname{mod} \Bbb Z).\tag3.10
$$
\endproclaim 

The argument in the proof is that since these constants $\tilde
a''_0(I)$ are locally 
determined, they depend, besides on $D$, only on the first $n$ homogeneous
terms in the symbols of the projections, and these are the same for
$\Pi $ and $\overline\Pi $.

The result of the theorem was shown in \cite{Woj} for the case where
$D=D_0$ in a selfadjoint product case, $\Pi =\Pi _+$ and $\Cal S$ is
of order $-\infty $, assuming that $D_{0,\Pi _+}$ is invertible. The
hypothesis on invertibility was removed by Lee in the appendix of
\cite{PW}; he shows moreover that $\zeta (D^2_{0,\Pi
_+},0)+\operatorname{dim}V_0(D_{0,\Pi _+} )=0$, so we conclude that 
$$
\zeta (D_{0,\Pi }^2,0)+\operatorname{dim}V_0(D_{0,\Pi } )=0,\text{
when }\Pi =\Pi _++\Cal S.\tag3.11$$

\medskip
\subhead 3.2 Results for eta functions \endsubhead

There are also such perturbation results for the eta function $\eta
(D_\Pi ,s)$, the meromorphic extension of $\Tr(D(D^*_\Pi D_\Pi
)^{-\frac{s+1}2})$, when (2.3) and
(3.4) hold:

\proclaim{Corollary 3.4} Assumptions of Corollary
{\rm 3.2}.
In
{\rm (2.24)} with $D_1=D$, 
one has that $\tilde a'_0(D)=a'_0(D)=0$, and 
$$\tilde a''_0(D)\;
(=a''_0(D)) \text{ is locally determined.}\tag3.12
$$ 
\endproclaim

In other words, the double pole
of $\eta (D_\Pi ,s)$ at $0$ vanishes and the residue at $0$ is
locally determined.
   
We underline that the hypotheses, besides  
(2.3), (3.4), only contain
requirements on principal symbols (namely the well-posedness of $\Pi
$ for $D$ and the
commutativity of the principal symbols of $\Pi $ and $A^2$). So the
result implies in particular that {\it the vanishing of the double pole of
the eta function is invariant under perturbations of $\Pi $ of order
$-1$} (respecting (3.4)). Earlier results have dealt with
perturbations of $\Pi _+$ of order $-\infty $ \cite{Woj}, or
perturbations of general $\Pi $ of order $-n$ \cite{G4}.

Now consider the simple pole of $\eta (D_\Pi ,s)$ at 0. Here we can
generalize the result of Wojciechowski 
\cite{Woj} on the 
regularity of the eta function after a
perturbation of order $-\infty $, to perturbations of order $-n$ of
general $\overline\Pi $:

\proclaim{Theorem 3.5} Assumptions of Theorem {\rm 3.3}.

In {\rm (2.24)} with $D_1=D$, 
the $\tilde a''_0(D)$-terms (and $a''_0(D)$-terms) for
$D_{\overline\Pi }^*D_{\overline\Pi }$ and 
$D_\Pi ^*D_\Pi $ are the same:$$
\tilde a''_0(D)(D_\Pi ^*D_\Pi )=\tilde a''_0(D)(D_{\overline\Pi }^*D_{\overline\Pi });\tag3.14
$$
in other words, $\operatorname{Res}_{s=0}\eta (D_\Pi
,s)=\operatorname{Res}_{s=0}\eta (D_{\overline\Pi },s)$.

In particular, if $\tilde a''_0(D)(D_{\overline\Pi }^*D_{\overline\Pi
})=0$ (this 
holds for $\Pi _+$ and for certain $\Pi (\theta )$ if $D$ equals
$D_0$ in a selfadjoint product case), then $\tilde a''_0(D)(D_\Pi
^*D_\Pi )=0$, i.e., 
the eta function 
$\eta (D_{\Pi },s)$ is regular at $0$.

\endproclaim 

 The argument is again that the local determinedness implies that
symbol changes below the first $n$ terms in the projection do not
enter in the constants.

The eta regularity for the case $\overline\Pi=\Pi _+$, $D$ equal to $D_0$ and
selfadjoint on $X$ with product volume element on $X_c$, was shown in
\cite{DW91} under 
the assumptions $n$ odd and $D$ compatible; this was extended to
general $n$ and not necessarily compatible $D$ in M\"uller \cite{M}. It
was shown for 
certain $\Pi (\theta )$ in \cite{BL, Th\. 3.12}.

The result on the regularity of the eta function at $s=0$ for
$(-n)$-order perturbations of the product case with $\overline \Pi
=\Pi _+$ has 
been obtained independently by Lei \cite{Le} at the same time as our
result, by another
analysis based on heat operator formulas. 

We refer to \cite{G6} for further discussions of $\tilde a''_0$. There
are some general results in \cite{G4} on the behavior of the other
coefficients under perturbations of  $\Pi $. 

\bigskip

\head 4. Perturbation of the interior operator \endhead

\medskip\subhead 4.1 General perturbation results \endsubhead

In this chapter, we discuss the behavior of
all the
logarithmic and nonlocal coefficients in (2.24) (not just the leading
ones) under perturbations of
$D$. In particular, it is interesting to compare with the special
situation of a
product-type operator
with $\Pi $ equal to $\Pi _>(A)$ plus a projection in the nullspace of
$A$. In that case, when $n$ is even, there are no logarithmic terms at
integer powers except possibly for $k=0$ if $\varphi \ne I$ (so they
occur only at half-integer powers
$(-\lambda )^{-k-\frac12}$); when $n$ is
odd, there are no logarithmic terms at all.  (This is known from
\cite{GS2}.) The results we now present are proved in \cite{G5}.

We have two kinds of results. One kind is a general statement in
the non-product case, that
when $D_1-D_2$ vanishes to a certain order at $\partial X$, then the
log-coefficients $a'_k$ up to a certain index are preserved when
$D_1$ is replaced by $D_2$, and the
$a''_k$-coefficients appearing together with them are perturbed only
by local terms. The arguments can
be used also when comparing the trace expansions
for a general $D$ with those for an associated $D_0$ of product type, under
suitable hypotheses on the volume form. 

The
other kind of result is concerned with perturbations of the product-case by
tangential operators commuting with $A$. Here it is found that in odd
dimensions, there is still a vanishing of all the
log-coefficients; on the other hand
nontrivial log-coefficients can be expected at both integer and
half-integer powers when $n$ is even.

We fix the boundary projection; it can be a general well-posed
projection in the
first kind of result, and in the second kind it is taken equal to
$\Pi _\ge(A)$. (Its perturbations follow the rules from Chapters 2
and 3, and from \cite{G4}.) 

Consider first the general non-product case.
When $D$ is given in the form (2.1) with (2.2) on $X_c$, we can
write, in the notation of \cite{G5},
$$
A_1=A+x_nP_1+P_0,\tag 4.1
$$
where $P_1$ is first-order tangential and $P_0$ is of order 0 and
{\it constant in }$x_n$ (since $A_{10}=A_{10}|_{x_n=0}+x_nA'_{10}$
where the last term may be absorbed in $x_nA_{11}$). The formal   
adjoint is 
$$
D^*=(-\partial _{x_n}+A'_1)\sigma ^*,\quad A'_1=A+x_nP^*_1+P'_0,\text{
on }X_c,\tag4.2
$$
where $P'_0=P^*_0-v^{-1}\partial _{x_n}v$; here $v(x)$ is the
function entering in the volume form $v(x)\,dx$. Conditions on the volume
form are needed when $D$ is compared with $D_0$ (coupled with the
volume form $v(x',0)\,dx'dx_{n}$) in Theorem 4.5 below.

Let $D_1$ and $D_2$ be two first-order
elliptic operators on
$X$ of non-product type (as in {\rm (2.1)} with {\rm (4.1)}), with
the same $\sigma $ and 
provided with the same well-posed 
boundary condition $\Pi \gamma 
_0u=0$. Let
$D_{1,\Pi }$ and $D_{2,\Pi }$ be the realizations defined by the
boundary condition $\Pi \gamma _0u=0$, and let $\Delta
_{i}=D_{i }^*D_{i }$, $\Delta
_{i,B}=D_{i,\Pi }^*D_{i,\Pi }$. 
Let $Q_{i,\lambda }$ be parametrices of the $\Delta _i$ on a
neighbouring manifold $\widetilde X$, and denote the resolvents of
the $\Delta _{i,B}$ by $R_{i,\lambda
}=(\Delta _{i,B}-\lambda )^{-1}$. 

We have the following general perturbation result:

\proclaim{Theorem 4.1}
Let $l$ be the largest nonnegative integer such that 
$$D_1-D_2=x_n^l\overline P_l\text{ on }X_c,\tag4.3
$$
for some tangential $x_n$-dependent first-order
differential operator $\overline P_l$. 
Let $F$ be a differential operator in
$E$ of 
order ${m'} $ and let $N>\frac {n+{m'} }2$. Consider $$
F(R_{2,\lambda }^N-R_{1,\lambda }^N)=(F(Q^N_{2,\lambda
}-Q^N_{1,\lambda }))_++
F\overline G^{(N)}_\lambda  .\tag4.4
$$

The $\psi $do part has an asymptotic trace expansion
$$\Tr [(F (Q^N_{2,\lambda }-Q^N_{1,\lambda }))_+]
\sim
\sum_{-n\le k<\infty  } \tilde p_{ k}(-\lambda )^{\frac{{m'} -k}2-N},
\tag4.5
$$
where $\tilde p_k=0$ for $k-m'+n$ odd.

The s.g.o.\ part has an asymptotic trace expansion
$$
\Tr [F\overline G^{(N)}_\lambda )]
\sim
\sum_{-n+1+l\le k< k_0 } \tilde g_{ k}(-\lambda )^{\frac{{m'}
-k}2-N}+ 
\sum_{k\ge k_0}\bigl({ \tilde g'_{ k}}\log (-\lambda
)+{\tilde g''_{ k}}\bigr)(-\lambda )^{\frac 
{{m'} -k}2-N},\tag4.6
$$
where $$
k_0=l+1 \text{ when $F$ is general},
k_0=m'+l+1 \text{ when $F$ is tangential on $X_c$.}
\tag4.7
$$
It follows that
$$
\Tr [F  (R^N_{2,\lambda }-R^N_{1,\lambda })]
\sim
\sum_{-n\le k< k_0 } \tilde c_{ k}(-\lambda )^{\frac{{m'}
-k}2-N}+ 
\sum_{k\ge k_0}\bigl({ \tilde c'_{ k}}\log (-\lambda
)+{\tilde c''_{ k}}\bigr)(-\lambda )^{\frac 
{{m'} -k}2-N},\tag4.8
$$
with $k_0$ as above. For $k\le l-n$, the $\tilde c_k$ vanish when
$k-m'+n$ is odd. 

The coefficients $\tilde c_k$ and $\tilde c'_k$ are locally determined.
 
\endproclaim 

The results carry over to similar results for the heat operators and
power operators associated with the $\Delta _{i,B}$. Alternatively,
we can formulate the results as follows:

\proclaim{Corollary 4.2}
Hypotheses and definitions as in Theorem
{\rm 4.1}. For the 
trace expansions$$
\aligned
\Tr (F R^N_{1,\lambda  })
&\sim
\sum_{-n\le k<0}  a_{ k}(-\lambda )^{\frac{{m'} -k}2-N}+ 
\sum_{k\ge 0}\bigl({ \tilde a'_{ k}}\log (-\lambda ) +{\tilde a''_{
k}}\bigr)(-\lambda ) ^{\frac{{m'} -k}2-N},\\
\Tr (F e^{-t\Delta _{1,B}})
&\sim
\sum_{-n\le k<0}  a_{ k}t^{\frac{k-m'}2}+ 
\sum_{k\ge 0}\bigl({ - a'_{ k}}\log t +{ a''_{
k}}\bigr)t ^{\frac{k-m'}2},\\
\endaligned$$
$$\aligned
\Tr (F \Delta _{1,B}^{-s})
&\sim
\sum_{-n\le k<0} \frac{ a_{ k}}{s+\frac{k-m'}2}-\frac{\Tr (F\Pi
_0(D_{1,\Pi }))}s\\
&\qquad+ 
\sum_{k\ge 0}\Bigl(\frac{  a'_{ k}}{(s+ \frac{k-m'}2)^2}+\frac{ a''_{
k}}{s+\frac{k-m'}2}\Bigr)
\endaligned
\tag4.9
$$ (the summation limit $0$ replaced by $m'$ if $F$ is tangential), the
 replacement of $D_1$ by
$D_2$ leaves the coefficients 
$\tilde a'_k$ and $a'_k$ invariant for $ k<k_0$. The other coefficients with
$k<k_0$ are modified only by local terms; those with
$ k\le l-n$ and $k-m'+n$ odd are invariant. 

\endproclaim

There are similar results for expansions associated with 
$D_{i,\Pi }R_{i,\lambda }$, $D_{i,\Pi }e^{-t\Delta _{i,B}}$
and $D_{i,\Pi }\Delta _{i,B}^{-s}$ (here the index $\Pi $ on the
factor in front can be omitted
since the resolvent and heat operator map into the domain):

\proclaim{Theorem 4.3} Hypotheses and definitions as in Theorem {\rm 4.1}. 
Let $\psi $ be a morphism from $E_1$ to $E$, and let $N>(n+m'+1)/2$. For the 
trace expansions$$
\aligned
\Tr (F \psi D_{1 }R^N_{1,\lambda  })
&\sim
\sum_{-n\le k<0} \tilde b_{ k}(-\lambda )^{\frac{{m'+1} -k}2-N}\\
&\qquad+ 
\sum_{k\ge 0}\bigl({ \tilde b'_{ k}}\log (-\lambda ) +{\tilde b''_{
k}}\bigr)(-\lambda ) ^{\frac{{m'+1} -k}2-N},\\
\Tr (F \psi D_{1 }e^{-t\Delta _{1,B}})
&\sim
\sum_{-n\le k<0}  b_{ k}t^{\frac{k-m'-1}2}+ 
\sum_{k\ge 0}\bigl({ - b'_{ k}}\log t +{ b''_{
k}}\bigr)t ^{\frac{k-m'-1}2},\\
\Tr (F \psi D_{1 }\Delta _{1,B}^{-s})
&\sim
\sum_{-n\le k<0} \frac{ b_{ k}}{s+\frac{k-m'-1}2}\\
&\qquad+ 
\sum_{k\ge 0}\Bigl(\frac{  b'_{ k}}{(s+ \frac{k-m'-1}2)^2}+\frac{ b''_{
k}}{s+\frac{k-m'-1}2}\Bigr)
\endaligned
\tag4.10
$$ (the summation limit $0$ replaced by $m'$ if $F$ is tangential), the
 replacement of $D_1$ by
$D_2$ leaves the coefficients 
$\tilde b'_k$ and $b'_k$ invariant for $ k<k_0$. The other coefficients with
$k<k_0$ are modified only by local terms; those with
$ k\le l-n$ and $k-m'+n$ even are invariant. 
\endproclaim 

Because of the factor $D_{i }$ in front, this is not a special
case of Theorem 4.1. Let us also mention, for the case $F=I$, the more
customary formulation of the third expansion, as in the last line of (2.24):

\proclaim{Corollary 4.4} Hypotheses and definitions as in Theorem {\rm 4.1}. 
Let $\psi $ be a morphism from $E_1$ to $E$, and let $k_0=l+1$. In
the eta function expansion  
$$\aligned
\eta (\psi,D_{1,\Pi },s)&=\Tr (\psi D_1 (D^*_{1,\Pi} D_{1,\Pi })^{-\frac{s+1}2})\\
&\sim
\frac 1{\Gamma (\frac{s+1}2)}\Bigl[\sum_{-n< k<0}\frac {2b_{k}}{s+ k}
+\sum_{k=0}^\infty
\Bigl(\frac {4b'_{k}}{(s+ k)^2}+\frac{2b''_{k}}{s+
k}\Bigr)\Bigr],
\endaligned \tag4.11
$$
a replacement of $D_1$ by
$D_2$ leaves the coefficients 
$\tilde b'_k$ and $b'_k$ invariant for $ k<k_0$. The other coefficients with
$k<k_0$ are modified only by local terms; those with
$ k\le l-n$ and $k+n$ even are invariant. 
\endproclaim 

The proofs are given in \cite{G5}; here we incorporate $D_{i,\Pi }$ and
${D_{i,\Pi }}^*$ in  larger skew-selfadjoint matrices$$
\Cal D_{i,\Cal B}=
\pmatrix 0&-{D_{i,\Pi }}^*\\ D_{i,\Pi }&0 
\endpmatrix  ,\tag4.12
$$
(as in \cite{GS1}) and study the difference of their resolvents,
using the calculus of 
\cite{G3} to handle the resulting singular Green operator term and
to find its expansion properties.

\medskip\subhead 4.2 Comparison with the product case \endsubhead

We can also compare the expansions for a given $D$ of non-product
type (2.1), (4.1), with the
expansions for an operator $D_0$ of product type having the form 
$D^0$ (2.5) on $X_c$. Here the volume form $v(x)\,dx$ for $D$ is
replaced by the volume form $v(x',0)\,dx$ for $D_0$ on $X_c$, so the
preceding results cannot immediately be applied. However, if for some
$l\ge 1$,
$$
P_0=0,\quad x_nP_1=x_n^lP_l, \quad 
\partial _{x_n}^jv(x',0)=0\text{ for }1\le j\le l,\tag4.13$$ 
then $D^*$ can be written in the form $$
D^*=(-\partial _{x_n}+A+x_n^lP'_l)\sigma ^*\text{ on }X_c;\tag4.14
$$
here $P_l$ and $P'_l$ are first-order tangential differential
operators.
Then the method of proof of the preceding results extends to show:

\proclaim{Theorem 4.5} Consider {\rm (4.9)} with $F$ equal to a
morphism $\varphi $, and {\rm (4.10)} with $F=I$ (so $m'=0$).

$1^\circ$ (The case $l=1$.) Assume that $P_0=0$ and $\partial
_{x_n}v(x',0)=0$. Then the 
coefficients $a'_0, a'_1$ (and $\tilde a'_0,\tilde a'_1$) in {\rm
(4.9)}  are the
same for the expansions defined for $D_{\Pi }$ and for
$D_{0,\Pi }$. The coefficients $a''_0, a''_1$ (and $\tilde
a''_0,\tilde a''_1$) differ in the two cases only by local terms.  

Moreover, in {\rm (4.10)}, the coefficients $b'_0, b'_1$ (and $\tilde
b'_0,\tilde b'_1$) are the same for $D_{\Pi }$ and for $D_{0,\Pi }$. 
The coefficients $b''_0, b''_1$ (and $\tilde
b''_0,\tilde b''_1$) differ in the two cases only by local terms.  

$2^\circ$ (The general case $l\ge 1$.) Assume that {\rm (4.13)} holds. Then in
{\rm (4.9), (4.10)}, the coefficients $a'_k$ and $b'_k$ for $0\le
k\le l$ (as well as $\tilde a'_k$ and $\tilde b'_k$ for $0\le k\le l$)
are preserved when $D_{\Pi }$ is replaced by $D_{0,\Pi }$. The
nonlocal coefficients behind them,  $a''_k, b''_k, \tilde a''_k,
\tilde b''_k$ with $0\le k\le l$ are only locally perturbed.

\endproclaim 
We also have the result that when $D-D_0=x_nP_1+P_0$ on $X_c$, the
zero-order operator $P_0$ not necessarily being 
0, then $a'_0$ is the same for $D_\Pi $ and $D_{0,\Pi }$, and $a''_0$
differs only by local terms. This was known from \cite{G1, GS1} in
cases where $\Pi $ equals $\Pi _\ge(A)$ or 
certain finite rank perturbations of it. However, $a'_1$ will in
general depend on 
$P_0$, as demonstrated in \cite {G5, Rem.\ 3.10}.

\medskip\subhead 4.3  Perturbation of the product case by commuting operators
\endsubhead 

The study of perturbations of the product case that commute with $A$
is somewhat different; here one can use functional calculus for the
operators near the boundary, expressing them as functions of $A$
(continuing the line of \cite{GS2}). The traces we study are reduced
to traces of pseudodifferential operators on the boundary, built up
of $A$ and its eigenprojections. Parity considerations play a great
role, because of the fact that $A$ and its integer powers have even-even
parity, whereas 
$|A|$ has even-odd parity, cf.\ (1.7)--(1.8).
By working out detailed formulas for the resolvent and its iterates
we were able to show in \cite{G5}:

\proclaim{Theorem 4.6} Assume that $D$ is a perturbation of $D^0$ as in
{\rm (4.1)} on $X_c$ such that the zero-order $x_n$-independent operator (morphism) $P_0$
commutes with $A$, and in the Taylor expansions on $X_c$,$$
x_nP_1(x_n)=\sum_{1\le k\le
K}x_n^kP_{1k}+x_n^{K+1}P'_{K+1}(x_n) \text{ for any }K, 
\tag 4.15$$
the tangential $x_n$-independent first-order differential operators
$P_{1k}$ commute with $A$. The product measure is used on $X_c$, and
$D$ is provided with the boundary condition$$
\Pi _\ge(A)\gamma _0u=0.\tag4.16
$$  
Let $F$ be a
differential operator in $E$ of order $m' $, and let $N>\frac{n+m' }2$.

{\bf If $n$ is odd,}  
the resolvent and heat operator, resp\. gamma times zeta function, associated with $\Delta _B$ have
trace expansions without logarithms, resp\. meromorphic extensions
without double poles:
$$\aligned
\Tr (F (\Delta _B-\lambda )^{-N})&\sim
\sum_{-n\le k< \infty } \tilde a_{ k}(-\lambda )^{\frac{m' -k}2-N}
,\\
\Tr (F  e^{-t\Delta _B})&\sim
\sum_{-n\le k<\infty } a_{ k}t^{\frac{k-m' }2} 
,\\
\Gamma (s)\zeta (F,\Delta _B,s)\equiv \Gamma (s)\Tr (F \Delta _B^{-s})&\sim
\sum_{-n\le k<\infty } \frac{a_{ k}}{s+\frac{k-m'}2}-\frac{\Tr(F \Pi
_0(\Delta _B))}{s} 
,
\endaligned\tag 4.17$$
where the coefficients are locally determined for $-n\le k<0$ (for
$-n\le k<m'$ if $F$ is tangential). Here $\tilde a_{-n}$ and $a_{-n}$
vanish if $m'$ is odd.
\endproclaim 

There are similar results for $\Tr (F\psi D (\Delta _B-\lambda
)^{-N})$ and its associated heat trace and power trace (an
eta-function), where $\psi $ is a morphism from $E_1$ to $E$.

In the case of a manifold of even dimension $n$, we first show that
for an $x_n$-independent and tangential differential operator $F$
taken together with a product case operator $D_0$, logarithmic terms
can appear at most at the power $(-\lambda )^{-N}$ and the
half-powers $(-\lambda )^{-N-k-\frac12}$, $k\in\Bbb N$; this is for
the expansion of $\Tr (F (\Delta _B-\lambda 
)^{-N})$, and there are corresponding statements for the other
trace expansions.

But now, when $D_0$ is replaced by a perturbation $D$ as in Theorem
4.6 (or when $F$ is $x_n$-dependent or non-tangential),
logs can appear at both integer and
half-integer powers in general.


\enddocument